\documentclass[a4paper,11pt]{article}
\usepackage{graphicx,color}
\usepackage{bmpsize,textcomp}
\usepackage{amsmath}
\usepackage{mathtools}
\usepackage{bbold}
\usepackage{bm}
\usepackage{fullpage}
\usepackage{url}
\usepackage{dsfont}
\usepackage{hyperref}

\title{\bf Computing first passage times \\ for Markov-modulated fluid models using\\ numerical PDE problem solvers}
\author{D. Bhaumik\footnote{Centrum Wiskunde \& Informatica, Science Park 123, 1098 XG Amsterdam, The Netherlands} , M.A.A. Boon\footnote{Eindhoven University of Technology, Department of Mathematics and Computer Science, P.O. Box 513, 5600 MB Eindhoven, The Netherlands.} , D. Crommelin\footnotemark[1] \footnote{Korteweg-de Vries Institute for Mathematics, University of Amsterdam, Science Park 105, 1098 XG Amsterdam, The Netherlands} , B. Koren\footnotemark[1] \footnotemark[2] , and B. Zwart\footnotemark[1] \footnotemark[2]}
\date{\today}

\graphicspath{{./},{./images/}}

\newcommand{\dd}{\text{d}}
\newcommand{\dt}{\text{d}t}
\newcommand{\dx}{\text{d}x}

\begin{document}

\maketitle
\title
\date

\begin{abstract}
A popular method to compute first-passage probabilities in continuous-time Markov chains is by numerically inverting their Laplace transforms.
Past decades, the scientific computing community has developed excellent numerical methods for solving problems governed by partial differential equations (PDEs), making the availability of a Laplace
transform not necessary here for computational purposes.
In this study we demonstrate that numerical PDE problem solvers are suitable for computing first passage times, and can be very efficient for this purpose. By doing extensive computational experiments, we show that modern PDE problem solvers can outperform numerical Laplace transform inversion,
even if a transform is available. When the Laplace transform is explicit (e.g. does not require the computation of an eigensystem), numerical transform inversion remains the primary method of choice.
\end{abstract}

\section{Introduction}
\par Much of the early research in applied probability and stochastic systems focused
on applying analytic techniques to queueing systems, deriving information about performance
measures in terms of Laplace transforms; the classic monograph on the single-server queue \cite{Cohen}
is a premier example. At that time, it was not clear how such transforms can be inverted, and other approaches,
in particular the usage of phase-type methods, gained traction after the seminal work \cite{Neuts}.
The applied relevance of Laplace transforms within applied probability got recognized in \cite{AbateWhitt},
where it was shown how Laplace transforms can be inverted numerically in a tractable manner.
A rigorous proof of convergence of the algorithms discussed in \cite{AbateWhitt} has recently
been established in \cite{Kuznetsov}.

\par The point we aim to make in this paper is that the field of applied probability and stochastic systems can utilize other computational methods.
In particular, we consider first-passage times of a class of structured Markov chains known as Markov-modulated fluid models.
Fluid queues are used to approximate discrete queueing models for dam theory \cite{moran1956probability}, telecommunication networks \cite{anick1982stochastic, elwalid1991analysis, mitra1987stochastic, schwartz1996broadband}, transportation systems \cite{newell2013applications}; to model forest fires \cite{stanford2005erlangized}, ruin probability \cite{asmussen2002erlangian, barges2013finite}, video streaming \cite{bosman2014optimal, kumar2007stochastic, wu2009queuing}, and more \cite{scheinhardt1998markov}. Steady state analysis of fluid queues is covered in detail in the literature, see e.g. \cite{asmussen1995stationary,  kulkarni1997fluid, rogers1994fluid}. Mean first passage times for fluid queues have been studied by \cite{asmussen1994sample, kulkarni2002mean}. The common approach to infer information about first-passage times
involves Laplace-Stieltjes transforms (LST's) \cite{boxma1998busy, field2010busy, narayanan1996first, ren1995transient}.
In \cite{narayanan1996first} it was proven that the cumulative distribution function of the first passage time of a Markov-modulated fluid queue model follows a \textit{hyperbolic partial differential equation} (PDE). In \cite{narayanan1996first}, the PDE problem is solved using LST's and eigenvalue theory.
The (numerical) inversion of the LST involves spectral decomposition of a matrix of size of the state space of the continuous-time Markov chain (CTMC), leading to cumbersome and expensive computations, thereby limiting the applicability of the method.

\par In this paper, we investigate alternative methods to numerically solve the PDE problem directly, without using the LST. Our focus is on cases with large cardinality of the state space of the CTMC.
Numerical integration of hyperbolic PDEs has been studied extensively in the scientific computing literature \cite{hundsdorfer2013numerical, leveque2002finite,smith1985numerical}.
The PDE central to this paper consists of a set of coupled linear PDEs with one space dimension. This set also arises in the study of chemically reacting species that are advected in a flow, and is therefore also known as an advection-reaction PDE \cite{hundsdorfer2013numerical, leveque2002finite}.
We propose to use first and second order accurate upwind space discretization schemes and different (explicit and implicit) time integration schemes. These schemes are well established in scientific computing but do not seem to be well-known in applied probability. We compare results from numerical integration of the PDE problem with Monte Carlo simulations to assess the accuracy. We find that the PDE schemes can efficiently and accurately compute the first passage time distribution and they can be applied to large state space systems. We also compare our PDE approach with the LST method proposed in \cite{narayanan1996first}.

\par In Section \ref{sec:system setup} we describe the system setup, buffer model, first passage time distribution and the PDE for the distribution function. In Section \ref{sec:PDE solvers} we propose different integration schemes for solving the PDE problem under consideration. In Section \ref{sec:LST method} we present an efficient implementation of the LST method proposed by \cite{narayanan1996first}. In Section \ref{sec:results} we compare the speed and accuracy of the different PDE integration schemes and the LST method, using different LST inversion algorithms. The paper finishes with some conclusions in Section \ref{sec:conclusions}.

\section{System setup}
\label{sec:system setup}
Let $\{Z(t), t \geq 0\}$ be a CTMC on a finite state space $\mathcal{S}$ with generator matrix $\bm{Q} = (q_{\alpha\beta})_{\alpha,\beta \in \mathcal{S}}$. The steady state probability distribution of $\{Z(t)\}$ is denoted by $\boldsymbol{\pi} = (\pi_1, \ldots, \pi_S)$, where $S = |\mathcal{S}|$. Let $r_\alpha$ be the net input rate with which the level of fluid in the buffer increases or decreases when the CTMC is in state $\alpha$. Let $B(t)$ be the fluid level in the buffer at time $t$. The dynamics of the fluid level is given by

\begin{align}
\label{eq: buffer model}
\frac{\dd B(t)}{\dd t} &=
\begin{cases}
r_{\alpha} ~~~~~~\qquad\qquad \text{if} ~ 0 \leq B(t) \leq ~  {B}^{\text{max}} ~ \text{and} ~ Z(t) = \alpha, \\
\text{max}(r_\alpha, 0) ~\qquad \text{if} ~ B(t)  = 0 ~ \text{and} ~ Z(t) = \alpha, \\
\text{min}(r_\alpha, 0) ~~\qquad \text{if} ~ B(t)  = {B}^{\text{max}} ~ \text{and} ~ Z(t) = \alpha. \\
\end{cases}
\end{align}

\noindent ${B}^{\text{max}}$ is the maximum fluid storage capacity of the buffer. Clearly, the fluid level increases when $r_\alpha > 0$ and
$B(t) < {B}^{\text{max}}$, and decreases when $r_\alpha < 0$ and $B(t) > 0$. Furthermore, $B(t) \in [0, {B}^{\text{max}}], ~\forall t$.
When $r_\alpha = 0$ the fluid level in the buffer remains unchanged.

\subsection{First passage time distribution}
\label{sec:first passage time}
We are interested in computing the distribution of the first passage (or hitting) time, i.e, the first time the buffer empties given that it started with some initial fluid level $x$ at time $t = 0$. In order to do so, we define a random variable $T$ which denotes the first time the buffer empties,

\begin{equation}
\label{eq:hitting time of battery empty}
T := \inf\{ t > 0 : B(t) = 0\}.
\end{equation}

\noindent The distribution function of $T$ given the initial fluid level in the buffer $x$, i.e., the probability that the buffer empties before time $t$ given that it started with initial fluid level $x$ is defined as

\begin{equation}
\label{eq:battery empty prob}
J(x,t) := \mathds{P}(T \leq t | B(0) = x).
\end{equation}

\par \noindent In the literature, slightly different variants of the above mentioned distribution function \eqref{eq:battery empty prob} are used depending on the initial conditions. A joint probability distribution of the buffer being empty before time $t$ and the state of the CTMC at time $T$ conditioned on the initial states of the buffer and the CTMC, is defined in \cite{narayanan1996first, ye2014computing}:

\begin{equation}
\label{eq:H}
H_{\alpha \beta}(x,t) := \mathds{P}(T \leq t, Z(T) = \beta | B(0) = x, Z(0) = \alpha),
\end{equation}

\noindent $\forall \alpha, \beta \in \mathcal{S}, x > 0$ and $t \geq 0$. In \cite{barbot2001distribution,jones2011fluid}, a probability distribution of $T$ conditioned on the initial state of the buffer fluid content and the CTMC is proposed:

\begin{equation}
\label{eq:K}
K_\alpha(x,t) := \mathds{P}(T \leq t| B(0) = x, Z(0) = \alpha) = \sum_{\beta = 1}^S H_{\alpha \beta}(x,t),
\end{equation}

\noindent $\forall \alpha \in \mathcal{S}, x > 0$ and $t \geq 0$. If we assume that the CTMC starts in stationarity, i.e. $\mathds{P}(Z(0) = \alpha) = \pi_\alpha, ~\forall \alpha \in \mathcal{S}$, we can compute the original probability distribution of $T$ given in  \eqref{eq:battery empty prob} from \eqref{eq:K} by

\begin{equation}
\label{eq:K to J}
J(x,t) = \sum_{\alpha = 1}^S \pi_\alpha K_\alpha(x,t).
\end{equation}

\noindent Hence, in this paper we focus on solving for $K_\alpha(x,t)$ which follows the following backward equation (see \cite{narayanan1996first} for the proof), $\forall \alpha \in \mathcal{S}$ and $\forall x, t > 0$,

\begin{equation}
\label{eq:K PDE}
\frac{\partial K_{\alpha}(x,t)}{\partial t} - r_\alpha \frac{\partial K_{\alpha}(x,t)}{\partial x} = \sum_{\beta = 1}^S q_{\alpha \beta} K_\beta(x,t).
\end{equation}

\noindent The initial and boundary conditions are given by:

\begin{itemize}				
\item $K_\alpha(x,0) = 0, \quad ~\forall x > 0, ~~\forall \alpha \in \mathcal{S}$.

\item $K_\alpha(0,0) = \begin{cases}
1, ~~~~\forall \{ \alpha : r_\alpha \leq 0 \}, \\
0, ~~~~\forall \{ \alpha : r_\alpha > 0\}.
\end{cases}$

\item $K_\alpha(0,t) =
1, ~~\forall \{ \alpha : r_\alpha \leq 0 \} \quad\text{and} \quad \forall t>0.
$

\item $K_\alpha(B^{\text{max}},t)$ follows an ordinary differential equation (ODE) at this boundary (when the buffer has reached $B^{\text{max}}$), i.e.,
\\
$\dfrac{\dd K_\alpha(B^{\text{max}},t)}{\dt} = \displaystyle\sum_{\beta = 1}^S q_{\alpha \beta} K_\beta(B^{\text{max}},t), ~~~~\forall \{ \alpha : r_\alpha > 0 \} \quad\text{and} \quad \forall t>0.$
\end{itemize}

Note that the individual PDEs in \eqref{eq:K PDE} can be coupled only through the summation term in the righthand side.

\section{Numerical schemes for solving PDEs}
\label{sec:PDE solvers}
In this section we summarize a few well-established PDE schemes for numerically integrating the PDE in \eqref{eq:K PDE}. Direct numerical integration for solving this PDE has not been deployed before in the field of queueing systems. The PDE in (\ref{eq:K PDE}) is a hyperbolic PDE, more specifically an \textit{advection-reaction} PDE with one space dimension (with the initial fluid level $x$ as space variable) \cite{hundsdorfer2013numerical, leveque2002finite}. 

\par For numerical integration of the PDE problem, the continuous space and time domains are discretized. The error of the  numerically integrated solution of the PDE depends on the discretization scheme used. In this section we discuss different space and time integration schemes. We present two \textit{upwind} space discretization schemes, namely the first order accurate upwind and the second order accurate upwind scheme. We present different schemes for time integrating the PDE problem, including one that is appropriate if the PDE is stiff
(in that case, the step size of explicit time integration methods must be kept very small to avoid numerical instability, eventually leading to excessively high computational cost).

\subsection{First order upwind scheme for spatial discretization}
\label{sec:Upwind scheme}

In this section we discuss a simple (first order accurate) space discretization scheme suitable for the PDE in (\ref{eq:K PDE}), the first order upwind finite difference scheme \cite{hundsdorfer2013numerical}.

\par
\begin{figure}[t!]
\hspace{0cm}
\includegraphics[scale = 0.75]{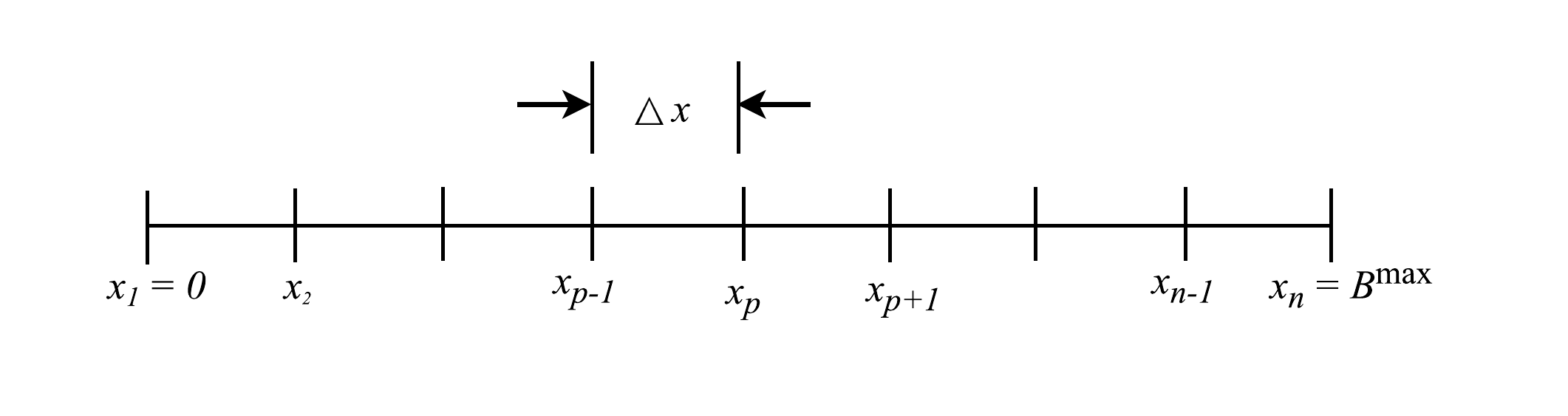}
\caption{Discretization of spatial domain for first order upwind scheme.}
\label{fig: upwind spacial discretization}
\end{figure}

 In our problem the space variable corresponds to the initial fluid level in the buffer. We discretize the interval $[0, B^{\text{max}}]$ using $n$ equidistant grid points with spacing $\Delta x = B^{\text{max}}/(n-1)$. The grid points are given by $\{ x_1 = 0, x_2, x_3, \ldots, x_n = B^{\text{max}} \}$ with $x_p = (p-1)\Delta x$ for $p = 1, \ldots, n$. The schematic representation of the grid point spacing and positions is shown in Figure \ref{fig: upwind spacial discretization}.

\noindent We discretize the space first, keeping the time still continuous. This will give us a set of coupled ordinary differential equations (ODEs). The first order upwind discretization scheme for the $p$-th space point is given by,

\begin{equation}
\label{eq:upwind space}
\frac{\dd K_\alpha^{p}}{\dt}  =
\begin{cases}
\frac{r_\alpha}{\Delta x} (K_\alpha^{p} - K_\alpha^{p-1}) + \displaystyle\sum_{\beta = 1}^S q_{\alpha\beta} K_\beta^{p} \qquad \forall \{\alpha : r_\alpha \leq 0 \},\\

\frac{r_\alpha}{\Delta x} (K_\alpha^{p+1} - K_\alpha^{p}) + \displaystyle\sum_{\beta = 1}^S q_{\alpha\beta} K_\beta^{p} \qquad \forall \{\alpha : r_\alpha > 0 \},
\end{cases}
\end{equation}

\noindent for $p = 2, \ldots, n-1$. In the above equation, $K_\alpha^p := K_\alpha(x_p)$ for $p = 1, \ldots, n$. The above equation can be written in a much more compact vectorial form. Let $\bm{K} = \left(K_1^{1}, \ldots, K_1^{n}, \ldots, K_S^{1},\ldots, K_S^{n}\right)^T$, then \eqref{eq:upwind space} can be written as
\begin{equation}
\label{eq:upwind space compact}
\frac{\dd \bm{K}}{\dt} = \bm{A} \bm{K},
\end{equation}
where $\bm{A}$ is a matrix of which the elements are determined by the space discretization and the boundary conditions. The structure of matrix $\bm{A}$ is given in Appendix \ref{App:first order upwind}. Note that \eqref{eq:upwind space compact} is a set of linear ODEs, which with proper initial conditions, can in principle be solved by

\begin{equation}
\label{eq:direct time integrate}
\bm{K}^{q+1} = e^{\bm{A} \Delta t} \bm{K}^{q}.
\end{equation}

\noindent In the above equation $\bm{K}^q := \bm{K}(t_q)$ and $t_{q+1} - t_{q} = \Delta t$.
In this paper we divide the time interval $[0,t_{\text{end}}]$ using $m$ time discretization points with spacing $\Delta t = t_{\text{end}}/(m-1)$. The $q$-th time point is given by $t_q = (q-1)\Delta t$ for $q = 1, \ldots, m$. At the initial time $t_1$, the solution, $K^1$, is known through the initial conditions.

\par Computing the exponential of the matrix $\bm{A}$ in \eqref{eq:direct time integrate} will be computationally very expensive as the size of the matrix is $(nS) \times (nS)$, hence it can be large. In the subsequent sections we describe efficient numerical time integration
schemes that can be used in conjunction with the first order upwind scheme.

\subsubsection{Time integration schemes}
\label{sec:time integration}
In this section we discuss two different time integration schemes, to be possibly applied to the first order upwind discretized set of ODEs. The first is the explicit fourth order \textit{Runge-Kutta} (RK4) scheme and the second is the implicit second order \textit{backward differentiation formula} (BDF2) scheme that is more suitable for stiff problems.

\paragraph{Explicit RK4 scheme:}
For numerical time integration, the integration time step $\Delta t$ needs to be specified. In case of a hyperbolic PDE discretized with a finite difference scheme, for numerical stability
of explicit time integration schemes, it is necessary that the \textit{Courant-Friedrichs-Lewy} (CFL) condition is satisfied \cite{hundsdorfer2013numerical, leveque2002finite}. For our problem the CFL condition yields
\begin{equation}
\label{eq:CFL}
\Delta t \leq \displaystyle\min_{\alpha\in \mathcal{S}} \left|\frac{\Delta x}{r_\alpha}\right|.   	
\end{equation}

\noindent
A more rigorous, but much more elaborate numerical stability requirement is:
\[
|\lambda_{p,\beta}|\,\Delta t\leq 1, \quad \forall p,\beta,
\]
where $\lambda_{p,\beta}$ are the eigenvalues of matrix $\bm{A}$, taking also into account the coupling term $\sum_{\beta=1}^S q_{\alpha\beta}\bm{K}_\beta^p$.

Given $\Delta x$ one must choose $\Delta t$ such that it satisfies the above condition. We deploy the explicit RK4 scheme:

\begin{equation}
\label{eq:RK4}
\bm{K}^{q+1} = \bm{K}^{q} +  \frac{\Delta t}{6} \left(\bm{p}_1 + 2\bm{p}_1 + 2\bm{p}_3 + \bm{p}_4\right),
\end{equation}

\noindent where
\begin{align*}
\bm{p}_1 &= \bm{A}\bm{K}^{q}, \\
\bm{p}_2 &= \bm{A}\left(\bm{K}^{q} + \frac{\Delta t}{2} \bm{p}_1\right),\\
\bm{p}_3 &= \bm{A}\left(\bm{K}^{q} + \frac{\Delta t}{2} \bm{p}_2\right),\\
\bm{p}_4 &= \bm{A}\left(\bm{K}^{q} +  \Delta t \bm{p}_3\right).
\end{align*}

\noindent The RK4 scheme is 4th order accurate, i.e.\ the local error is	 ${\mathcal{O}(\Delta t^{4})}$. The computational complexity of the RK4 scheme along with the first order upwind scheme is $\mathcal{O}(mnS^2)$. Note that $\bm{A}$
is a sparse matrix with $\mathcal{O}(nS^2)$ nonzero elements. The combination of RK4 and first-order upwind is not optimal from the perspective of accuracy, because it is unbalanced; fourth order accuracy versus first order accuracy only. From the perspective of numerical stability though, it is a good combination.

\par The computational cost of an explicit scheme such as RK4 can become high in case of stiffness. In that case, using an implicit time integration scheme is beneficial. A system of ODEs is said to be stiff if the time step length for explicit schemes becomes very restrictive (i.e., very small) due to stability requirements
even though the solution varies only slowly in time.
For the system under consideration here, stiffness can occur if some of the rates $r_\alpha$ (or $q_{\alpha \beta}$) have different orders of magnitude. To overcome this problem we resort to the implicit time integration scheme discussed in the following paragraph.

\paragraph{Implicit BDF2 scheme:} The backward differentiation formula (BDF) is a family of implicit schemes for integrating ordinary differential equations. The $k$-th order BDF scheme uses $k$ previous time step approximations of the solution for computing the present approximation of the solution. We use the second order BDF (BDF2) scheme which is given by \cite{suli2003introduction}

\begin{equation}
\label{eq:BDF2}
\bm{K}^{q+1} = \frac{4}{3}\bm{K}^{q} - \frac{1}{3}\bm{K}^{q-1} + \frac{2}{3} \Delta t \,\bm{A} \bm{K}^{q+1}.
\end{equation}

\noindent The method is implicit, and to compute $\bm{K}^{q+1}$ given $\bm{K}^{q}, \bm{K}^{q-1}$ a system of linear equations needs to be solved. This makes the method computationally expensive. However, the method is less restrictive on $\Delta t$ for stability \cite{suli2003introduction}. In order to solve the system of linear equations in \eqref{eq:BDF2} we use the LU decomposition method (worst case complexity $\mathcal{O}((nS)^3)$ \cite{gilbert1988sparse}) implemented in MATLAB. In practice, the sparsity of the matrix will reduce the computational complexity of the LU decomposition. The degree of reduction is hard to assess because it depends on the sparsity pattern. The complexity of BDF2 with first order upwind scheme becomes $\mathcal{O}(mn^2S^2 + (nS)^3)$, where $m$ is the number of time steps and $n$ the number of space discretization points.

\subsection{Higher order upwind scheme using flux limiters}
\label{sec:2nd order upwind with flux}
In Section \ref{sec:Upwind scheme} we described the method of solving the PDE \eqref{eq:K PDE} using the first order accurate finite difference upwind space discretization scheme. The drawback of the this scheme is that the error in the solution is $\mathcal{O}(\Delta x)$. In order to achieve higher order accuracy in the numerically integrated solution of the PDE problem one can use higher order schemes. Hence, in this section we discuss among others the second order accurate upwind scheme for space discretization which has a solution error $\mathcal{O}(\Delta x^2)$.

\par The second order upwind scheme may lead to oscillations around discontinuities in the solution domain, though. In our problem the distribution function in \eqref{eq:K} has jumps, leading to discontinuities in the solution. These jumps occur because there is a positive probability that the CTMC starts and remains in a state with negative rate ($r_\alpha < 0$) for any finite amount of time and hence there is a positive probability that the buffer hits zero exactly at time $|x r_\alpha|$, where $x$ is the initial fluid level in the buffer. To limit the values of the spatial derivatives to realistic values around discontinuities of the solution, the concept of \textit{flux limiters} is used.

\subsubsection{Second order upwind}
\label{sec:2nd upwind scheme}
\begin{figure}[t!]
\hspace{0cm}
\includegraphics[scale = 0.75]{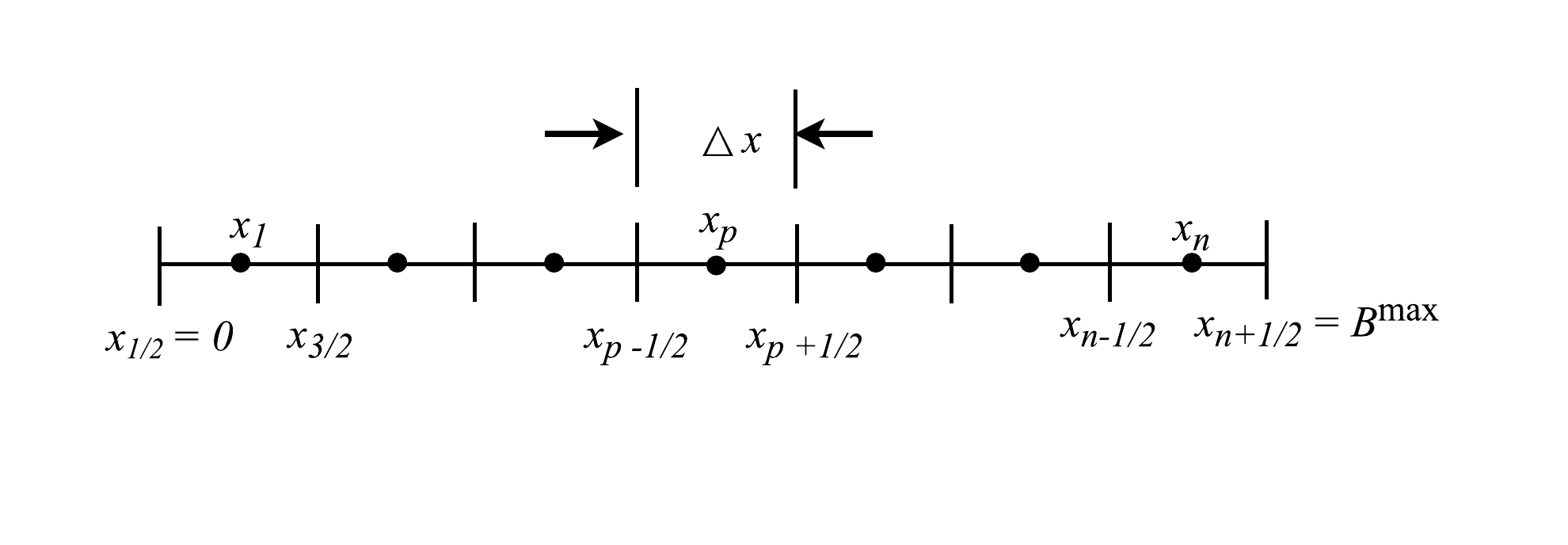}
\caption{Cell-centered discretization of spacial domain for second order upwind scheme.}
\label{fig: second upwind spacial discretization}
\end{figure}

The PDE in \eqref{eq:K PDE} is spatially approximated on cell-centered grid points $x_p$ (see Figure \ref{fig: second upwind spacial discretization}) by the semi-discrete conservation form \cite{hundsdorfer1995positive}

\begin{equation}
\label{eq: K second order}
\frac{\dd K_\alpha^p(t)}{\dt} - \frac{1}{\Delta x} \left(F_
 \alpha^{p + 1/2}(t) - F_\alpha^{p - 1/2}(t)\right) = \sum_{\beta = 1}^S q_{\alpha\beta} K_\beta^p(t), ~~ \forall \alpha \in \mathcal{S}.
\end{equation}

\noindent In the above equation, $K_\alpha^p(t)$ is the continuous-time approximation of $K_\alpha(x_p,t)$ at the grid point $x_p = (p - \frac{1}{2})\Delta x$ for $p = 1, \ldots, n$, where $n$ is the number of cell centered grid points. $F_\alpha^{p+1/2}$ is the numerical flux at the right vertex of the $p^{th}$ cell and it is the approximation of the analytical flux at cell vertex $x_{p+ 1/2} = \frac{1}{2} (x_{p+1} + x_p)$ for $p = 1, \ldots, n-1$. The boundary vertices are $x_{1/2} = 0$ and $x_{n + 1/2} = B^{\text{max}}$.

\par For the second order upwind scheme the numerical flux is given by \cite{hundsdorfer1995positive}

\begin{equation}
\label{eq:flux}
F_\alpha^{p+1/2} = \begin{cases}
r_\alpha\left(\frac{3}{2} K_\alpha^p - \frac{1}{2} K_\alpha^{p-1}\right) \qquad\forall \{ \alpha : r_\alpha \leq 0 \}, \\
r_\alpha\left(\frac{3}{2} K_\alpha^{p+1} - \frac{1}{2} K_\alpha^{p+2}\right) \quad\forall \{ \alpha : r_\alpha > 0 \}.
\end{cases}
\end{equation}
\noindent For other higher order (for example third and fourth order) discretization schemes, see \cite{hundsdorfer1995positive}.

\subsubsection{Positivity of the second order upwind scheme}
\label{sec:positivity}
As mentioned earlier, the second order upwind scheme may suffer from oscillations around discontinuities in the solution domain. It implies that the scheme can lead to negative values in the numerical solution. This is unwanted as we are approximating a distribution function whose values must lie in $[0,1]$. In this subsection we give an explanation for this behavior of the scheme. We follow \cite{hundsdorfer1995positive} closely in this section.

\par The positivity rule says that for any non-negative initial solution $( K_\alpha^p(t_0) \geq 0, ~\forall p )$ the solution $K_\alpha^p(t)$ which evolves in time remains non-negative, $\forall t \geq t_0$. A scheme is positive if and only if, for $K_\alpha^p(t)=0$, $\forall p$ and $\forall t \geq t_0$ \cite{hundsdorfer1995positive}:
\begin{equation}
\label{eq:positivity}
K_\alpha^{p'}(t) \geq 0, ~~\forall p' \neq p \implies \frac{\dd K_\alpha^p(t)}{\dt} \geq 0.
\end{equation}

\noindent Equation \eqref{eq:positivity} means that if the solution in any cell $p' \neq p$ increases, then the solution in cell $p$ should also increase, or remain the same. It should be noted that the first order upwind scheme (see Section \ref{sec:Upwind scheme}) satisfies this rule. For the second order upwind scheme we plug in the values of the numerical fluxes from \eqref{eq:flux} into \eqref{eq: K second order} and we get $\forall \{ \alpha : r_\alpha < 0 \}$
\begin{equation}
\frac{\dd K_\alpha^p(t)}{\dt} = \frac{|r_\alpha|}{2 \Delta x} \left(-K_\alpha^{p-2}(t) + 4 K_\alpha^{p-1}(t) -3 K_\alpha^{p}(t)\right) + \sum_{\beta = 1}^S q_{\alpha\beta} K_\beta^p(t).
\end{equation}

We see that the coefficient of $K_\alpha^{p-2}(t)$ is negative and hence when $K_\alpha^{p-2}(t)$ increases $\dfrac{\dd K_\alpha^{p}(t)}{\dt}$ decreases. Therefore, the positivity rule is not satisfied by the second order upwind scheme. The $q$-terms in the above equation pose no problems for the positivity rule because of the structure of the $\bm{Q}$ matrix, i.e, the non-diagonal terms of the matrix are non-negative. Analogous arguments show that the positivity rule is also not satisfied for $\frac{\dd K_\alpha^p(t)}{\dt} ~~ \forall \{ \alpha : r_\alpha \geq 0 \}$.

\subsubsection{Flux limiters}
\label{sec:flux limiters}
We have seen in the previous section that the second order upwind scheme does not preserve non-negativity of the time evolving solution. To ensure non-negativity of the solution we apply the \textit{flux limiter} method. The expression for the flux with the limiter is given by

\begin{equation}
\label{eq:flux limiter}
F_\alpha^{p+1/2} = \begin{cases}
r_\alpha\left(K_\alpha^p + \frac{1}{2}\phi_{p+1/2}(K_\alpha^p - K_\alpha^{p-1})\right) \hspace{3mm}\qquad\forall \{ \alpha : r_\alpha \leq 0 \}, \\
r_\alpha\left(K_\alpha^{p+1} + \frac{1}{2}\phi_{p+1/2}(K_\alpha^{p+1} - K_\alpha^{p+2})\right) \quad\forall \{ \alpha : r_\alpha > 0 \},
\end{cases}
\end{equation}

\noindent where $\phi_{p+1/2} := \phi(f_{p+1/2})$. Here $\phi$ is a non-linear function called the \textit{limiter function}, and $f_{p+1/2}$ is the ratio of consecutive solution gradients. Both $\phi$ and $f$ will be specified later in this section. The limiter  is applied such that the solution remains non-oscillatory in the discontinuous part of the solution and in the smooth part the higher order scheme is applied \cite{hundsdorfer1995positive,leveque2002finite}. Note that with $\phi = 0$, we get the first order upwind scheme. We show that the flux limited expression in \eqref{eq:flux limiter} preserves positivity of the semi-discrete scheme in \eqref{eq: K second order}. For $\{ \alpha : r_\alpha < 0 \}$  we get from \eqref{eq:flux limiter} and \eqref{eq: K second order}

\begin{align}
\label{eq:flux prove1}
\frac{\dd K_\alpha^p}{\dt} &= -\frac{|r_\alpha|}{\Delta x} \left((K_\alpha^p - K_\alpha^{p-1}) + \frac{1}{2} \phi_{p + 1/2} (K_\alpha^p - K_\alpha^{p-1}) - \frac{1}{2} \phi_{p - 1/2} (K_\alpha^{p-1} - K_\alpha^{p-2})\right) \nonumber
\\
& + \sum_{\beta = 1}^S q_{\alpha\beta} K_\beta^p.
\end{align}

\noindent Let
\begin{equation}\label{eqn:fp}
f_{p-1/2} = \dfrac{K_\alpha^p - K_\alpha^{p-1}}{K_\alpha^{p-1} - K_\alpha^{p-2}}.
\end{equation}
Note that with the choice \eqref{eqn:fp}, and hence also $f_{p+1/2} = \dfrac{K_\alpha^{p+1} - K_\alpha^{p}}{K_\alpha^{p} - K_\alpha^{p-1}}$, the space discretization is no longer fully one-sided upwind pers\'e, but may also be upwind biased (partly second order upwind, partly second order central), enabling besides positivity also third order accuracy.
Assuming $f_{p-1/2} \neq 0 $, i.e., $K_\alpha^p - K_\alpha^{p-1} \neq 0$, we get from \eqref{eq:flux prove1}:

\begin{equation}
\label{eq:flux prove2}
\frac{\dd K_\alpha^p}{\dt} = -\frac{|r_\alpha|}{\Delta x} \left( (1 + \frac{1}{2} \phi_{p+1/2}) - \frac{1}{2} \frac{\phi_{p-1/2}}{f_{p-1/2}} \right) (K_\alpha^p - K_\alpha^{p-1}) + \sum_{\beta = 1}^S q_{\alpha\beta} K_\beta^p.
\end{equation}

\noindent Applying the positivity rule from \eqref{eq:positivity} to \eqref{eq:flux prove2}, we find that the term in the square bracket in the above equation should be non-negative and this leads to

\begin{equation}
\label{eq:flux prove 3}
 \frac{\phi_{p-1/2}}{f_{p-1/2}} -\phi_{p+1/2} \leq 2.
\end{equation}

\noindent If $f_{p-1/2} = 0$, i.e., $K_\alpha^p - K_\alpha^{p-1} = 0$, then \eqref{eq:flux prove1} leads to $\frac{\dd K_\alpha^p(t)}{\dt} > 0 $ if $\phi_{p-1/2} = 0$. Hence, we assume that $\phi_{p-1/2} = 0$ if $f_{p-1/2} \leq 0$. Also assuming $0 \leq\phi_{p-1/2}, \phi_{p+1/2} \leq \delta$, for any constant $\delta > 0$, then \eqref{eq:flux prove 3} is satisfied if $\phi_{p - 1/2} \leq 2 f_{p-1/2}$. Hence, $\phi_{p-1/2}$ and $\phi_{p+1/2}$ can be expressed as functions of $f_{p-1/2}$ and $f_{p+1/2}$. Thus we have,

\begin{equation}
\label{eq:tvd}
\begin{cases}
\phi(f) = 0  ~\quad \text{for } f \leq 0, \\
\phi(f) \leq 2f  \quad \text{for } f > 0, ~~\text{and} \\
\phi(f) \leq 2    ~\quad \text{for } f > 0, ~
\end{cases}
\end{equation}
where we have chosen $\delta=2$.
\noindent From the conditions in \eqref{eq:tvd} the limiter must be within or at the boundary of the so-called Sweby TVD (total variation diminishing) domain (see Fig.\ 1a of \cite{sweby1984high}). One gets the same results $\forall \{ \alpha : r_\alpha \geq 0 \}$.

\par A wide range of limiters that fall in this TVD region has been proposed in literature, see \cite{waterson2007design, alhumaizi2007flux} for an overview and comparisons of different limiters. For solving our problem we use the Koren limiter \cite{koren1993robust}, given by

\begin{equation}
\label{eq:barry limiter}
\phi(f) = \text{max}(0, \text{min}(2f, \frac{1+2f}{3}, 2)).
\end{equation}
An advantage of limiter \eqref{eq:barry limiter} is that it makes the space discretization third-order accurate.

\subsubsection{Time integration scheme}
\label{sec:RK3b}
For time integration we cannot use the RK4 scheme (see Section \ref{sec:time integration}) with the flux limited scheme. This is because for the RK4 scheme to satisfy the TVD conditions, the step size must be zero (see equation (33) of
 \cite{hundsdorfer1995positive}). An advantage of the use of the RK3b scheme and the space discretization with limiter \eqref{eq:barry limiter} is that the orders of accuracy of both methods are the same; both third order. Hence, we use the explicit third order accurate Runge-Kutta 3b (RK3b) scheme described and suggested in \cite{hundsdorfer1995positive}. We describe the scheme in this section. Let the space discretized $K_\alpha(x)$ be given by $K_\alpha^{p} := K_\alpha(x_p)$, where $p$ is the index of the space discretization point (see Section \ref{sec:2nd upwind scheme}). Let $\bm{K} = [K_1^{1}, \ldots, K_1^{n}, \ldots,
K_S^{1}, \ldots, K_S^{n}]^T$. Then we can rewrite
the system of equations \eqref{eq: K second order} by expanding $F_
\alpha^{p+1/2}$ using \eqref{eq:flux limiter} in the form
\begin{equation}
\label{eq:second upwind space compact}
\frac{\dd \bm{K}}{\dt} = \bm{A}(\bm{K}) \bm{K},
\end{equation}

\noindent where the matrix $\bm{A}(\bm{K})$ is a function of $\bm{K}$ because of the flux-limiter (which is a function of $\bm{K}$). As can be seen in \eqref{eq: K second order} and \eqref{eq:flux limiter}, the elements of $\bm{A}(\bm{K})$ are determined by the third-order accurate space discretization, the boundary conditions and the source terms. For details of the space discretization see Appendix \ref{App:second order upwind}.
Note that unlike the upwind scheme described in Section \ref{sec:Upwind scheme} and \eqref{eq:upwind space compact} the matrix $\bm{A}$ now has to be evaluated at each of the three stages of the RK3b scheme. 
The explicit RK3b is given by
\begin{equation}
\label{eq:RK3b}
\bm{K}^{q+1} = \bm{K}^{q} +  \frac{\Delta t}{6} \left(\bm{p}_1 + \bm{p}_2 + 4\bm{p}_3\right),
\end{equation}

where
\begin{align*}
\bm{p}_1 &=  \bm{A}(\bm{K}^q) \bm{K}^q,\\
\bm{p}_2 &= \bm{A}(\bm{K}^{q} + \Delta t \bm{p}_1)\left(\bm{K}^{q} + \Delta t \bm{p}_1\right),\\
\bm{p}_3 &=  \bm{A}\left(\bm{K}^{q} + \frac{\Delta t}{4} \left(\bm{p}_1 + \bm{p}_2\right)\right)\left(\bm{K}^{q} + \frac{\Delta t}{4} \left(\bm{p}_1 + \bm{p}_2\right)\right),
\end{align*}

\noindent where $q$ is the discretized time index and the time step $\Delta t$ is chosen in such a way that the positivity requirement is satisfied, which latter is more restrictive than the numerical stability requirement. The positivity requirement on $\Delta t$, formula (30) in \cite{hundsdorfer1995positive}, reads:
\begin{equation}\label{cfl}
\Delta t = \leq \frac{1}{1+\delta/2}\displaystyle\min_{\alpha \in \mathcal{S}} \left|\frac{\Delta x}{r_\alpha}\right|.
\end{equation}
With $\delta=2$, this yields $\Delta t\leq \frac12 \displaystyle\min_{\alpha \in \mathcal{S}} \left|\frac{\Delta x}{r_\alpha}\right|$. For safety -- the equations have become nonlinear due to the limiter -- we take $\Delta t = \frac14  \displaystyle\min_{\alpha \in \mathcal{S}} \left|\frac{\Delta x}{r_\alpha}\right|$.  The computational complexity of the scheme is $\mathcal{O}(nmS^2)$, where as before, $n$ is the number of space discretization points, $m$ the number of time discretization points and $S$ the cardinality of the state space of the CTMC.

\section{Laplace-Stieltjes transform method}
\label{sec:LST method}
In queueing theory, the common approach to compute first passage times for fluid queues is by the Laplace-Stieltjes transform (LST) \cite{boxma1998busy, field2010busy, narayanan1996first, ren1995transient}. In this section, we start with a brief discussion of the LST method presented in \cite{narayanan1996first} for calculating \eqref{eq:K}. To overcome some numerical challenges, we develop an efficient algorithm to compute and invert the LST at the end of this section.

\subsection{General idea}

The starting point in \cite{narayanan1996first} is the distribution function defined in \eqref{eq:H}. The PDE satisfied by $H_{\alpha \beta}(x,t)$, $\forall \{ \beta : r_\beta < 0 \}, \forall \alpha \in \mathcal{S}$ and $\forall x, t > 0$ is

\begin{equation}
\label{eq:PDE}
\frac{\partial H_{\alpha\beta}(x,t)}{\partial t} - r_\alpha \frac{\partial H_{\alpha\beta}(x,t)}{\partial x} = \sum_{\gamma = 1}^S q_{\gamma \beta} H_{\gamma\beta}(x,t),
\end{equation}
with appropriate boundary and initial conditions discussed in \cite{narayanan1996first}. Note that $H_{\alpha\beta}(x,t) = 0$, $\forall \{ \beta : r_\beta > 0 \}$, and $ \forall \alpha \in \mathcal{S}$. Taking the LST,
$\tilde{H}_{\alpha\beta}(x,w) := \int_{0}^\infty e^{-wt} dH_{\alpha\beta}(x,t)$,
the vector $\tilde{\bm{H}}_\beta = [\tilde{H}_{1\beta}, \ldots \tilde{H}_{S\beta}]$ satisfies the following ODE

\begin{equation}
\label{eq:ODE}
\bm{R} \frac{\dd \tilde{\bm{H}}_\beta(x,w)}{\dx} = (w\bm{I} - \bm{Q}) \tilde{\bm{H}}_\beta(x,w),
\end{equation}
where $\bm{R} = \text{diag}(r_1, \ldots r_S)$. The following functional form is assumed
\begin{equation}\label{eq28}
\tilde{\bm{H}}_\beta(x,w) = e^{s(w) x} \bm{\Phi}(w).
\end{equation}
Substituting \eqref{eq28} in \eqref{eq:ODE} we obtain
\begin{equation}
\label{eq:eig1}
\bm{R} s(w)\bm{\Phi}(w) = (w \bm{I} - \bm{Q}) \bm{\Phi}(w),
\end{equation}

\noindent where $s(w)$ is a scalar and $\bm{\Phi}(w)$ a vector that are both to be determined. The above equation can be re-written as

\begin{equation}
\label{eq:eig2}
(\bm{Q} + s(w)\bm{R} -w\bm{I})\bm{\Phi}(w) = 0.
\end{equation}

\noindent Hence, the problem boils down to solving an eigenvalue problem, i.e., we need to find the roots of $\Delta(s,w) = \text{det}(\bm{Q} + s(w)\bm{R} -w\bm{I}) = 0$. Assuming that the diagonal elements of $\bm{R}$ are nonzero, let $s_k(w)$ for $k = 1, \ldots, S$ be the roots of $\Delta(s,w)$ and $\bm{\Phi}^k(w)$ the corresponding eigenvectors. We have

\begin{equation}
\label{eq:til_H final}
\tilde{H}_{\alpha\beta}(x,w) = \sum_{k = 1}^S a_{k \beta} e^{s_k(w)x} \phi^k_\alpha(w),
\end{equation}

\noindent where $\phi^k_\alpha(w)$ are the elements of $\bm{\Phi}^k(w)$. The coefficients $a_{k \beta}$ are obtained from the following initial and boundary conditions $\forall \{ \beta : r_\beta < 0 \}$:

\begin{itemize}

\item $\tilde{H}_{\alpha\beta}(0,w) = \begin{cases}
1 ~~~~\text{if} \quad \alpha = \beta,\\
0 ~~~~ \text{if} \quad \alpha \neq \beta \quad \text{and} \quad r_\alpha < 0.
\end{cases}$

\item $\tilde{H}_{\alpha\beta}(B^{\text{max}},w) =  \sum\limits_{m \neq \alpha} \frac{q_{\alpha m}}{-q_{\alpha\alpha} + w} \tilde{H}_{m \beta}(B^{\text{max}},w), ~~ \forall \{ \alpha : r_\alpha > 0 \}$.
\end{itemize}

\noindent From the above boundary conditions we can obtain the full expression for $\tilde{H}_{\alpha\beta}(x,w)$ from which we can obtain

\begin{equation}
\label{eq:til_K}
\tilde{K}{_\alpha}(x,w) = \sum\limits_{\beta = 1}^S \tilde{H}_{\alpha\beta}(x,w), ~~ \forall \alpha \in \mathcal{S}.
\end{equation}

\noindent Finally, to obtain $K_\alpha(x,t)$ as in \eqref{eq:K}, we need to numerically invert the transformed solution $\tilde{K}{_\alpha}(x,w)$. An advantage of the LST method is that the moments of the first passage time can also be obtained very efficiently.

\subsection{Numerical algorithm}

Summarizing the LST method, one first needs to solve the eigenvalue problem in (\ref{eq:eig2}). Then one solves the system of linear equations
\eqref{eq:til_H final} to obtain the coefficients $a_{k \beta}$. Finally, for applying the result to practical problems, one must numerically invert the LST to obtain the solution $K_\alpha(x,t)$ in the time domain.

The first decision that one has to make, when searching for the most efficient implementation, is whether to solve the eigenvalue problem analytically or numerically. The disadvantage of the analytical approach is that it requires a long computing time to find closed-form expressions for the eigenvalues. However, this procedure has to be executed only once, and with the resulting expressions the LST inversion becomes much faster. Still, the expressions tend to grow in size exponentially and even for moderately low values of $S$ it rapidly becomes infeasible to determine them. The alternative is to numerically determine the eigenvalues. This needs to be done for every value of $w$ where the LST is to be evaluated. In any case, much time is gained by rewriting Equations \eqref{eq:til_K} and \eqref{eq:til_H final} in matrix form. For notational reasons, we assume without loss of generality that the rates $r_\alpha$ are ordered from largest rate to smallest rate. Let $S_+$ denote the number of nonnegative rates, and $S_-$ the number of negative rates.

Since we only need the sum (over $\beta$) of $\tilde{H}_{\alpha\beta}(x,w)$ to determine $\tilde{K}{_\alpha}(x,w)$, we introduce the $S$-dimensional vector $\bm{\tilde{a}}(w)$ with elements $\tilde{a}_\alpha=\sum_{\beta=1}^S a_{\alpha\beta}$. Note that $\tilde{a}_\alpha=\sum_{\beta=1+S_+ }^S a_{\alpha\beta}$, since $H_{\alpha\beta}=a_{\alpha\beta}=0$ for $\beta\leq S_+$.

We rewrite \eqref{eq:til_K} as
\begin{equation}
\label{eq:til_K2}
\tilde{K}{_\alpha}(x,w) = \sum_{k = 1}^S e^{s_k(w)x} \phi^k_\alpha(w) \tilde a_{k} ,
\end{equation}
where the $\tilde{a}_{k}$ follow from the following two sets of equations. For $\alpha \leq S_+$,
\begin{equation}
\label{eq:set1}
\sum_{m=1}^S M_{\alpha,m}(w) \sum_{k=1}^S \tilde{a}_k e^{s_k(w)B^\text{max}} \phi^k_\alpha(w)=0,
\end{equation}
with $\bm{M}(w)$ an $S_+ \times S$ dimensional matrix with elements
\[
M_{\alpha,m}(w) =\begin{cases}
1&\quad \alpha=m,\\
\frac{q_{\alpha m}}{q_{\alpha\alpha}-w}&\quad \alpha\neq m.
\end{cases}
\]
For $\alpha > S_+$,
\begin{equation}
\label{eq:set2}
\sum_{k=1}^S \tilde{a}_k \phi^k_\alpha(w) = 1.
\end{equation}

Denote by $\bm{\Phi}(w)$ the matrix of eigenvectors (where $\bm{\Phi}^k(w)$ constitutes the $k$-th column), and by $\bm{\Phi_-}(w)$ the last $S_-$ rows of $\bm{\Phi}^k(w)$ (corresponding to the negative rates). 
The system of equations \eqref{eq:set1}--\eqref{eq:set2} can now be written in matrix form, as follows
\begin{equation}
\label{eq:systemofeqnsmatrixform}
\left(\begin{array}{c}
\\
\bm{M}(w)\bm{\Phi}(w) \text{diag}(e^{s_1(w)x}, \dots, e^{s_S(w)x})\\
\\
\hline
\\
\bm{\Phi_-}(w)
\\
\end{array}\right)\bm{\tilde{a}}(w) =
\left(\begin{array}{c}
0\\
\vdots\\
0\\
\hline
1\\
\vdots\\
1
\end{array}\right).
\end{equation}
This can be solved numerically for $\bm{\tilde{a}}(w)$ extremely fast and efficiently. The LST $\tilde{K}{_\alpha}(x,w)$ can be obtained from \eqref{eq:til_K2} or, more efficiently, from
\[
\tilde{K}{_\alpha}(x,w) = \bm{\Phi}(w) \text{diag}(e^{s_1(w)x}, \dots, e^{s_S(w)x})\bm{\tilde{a}}(w).
\]

\section{Results}
\label{sec:results}
In this section we compare the different methods to compute the distribution of the first passage times, i.e, the first time the buffer empties given that it started with some initial fluid level $x$ at time $t = 0$. We compare the accuracy of the different methods and discuss possible numerical issues encountered while running various examples.
Although we also discuss the computing times for the different algorithms, it is very difficult to compare them since they depend strongly on the chosen parameters.
Furthermore, the PDE methods compute $J(x,t)$ for all space discretization points ($x_1, ..., x_n$) at once, whereas the LST method gives the result at a single point in space.
In all examples, we determine $J(x,t)$ for a single $x$ value only. If one is interested in obtaining the distribution of the first passage time for multiple values of $x$, the PDE methods are (almost) always preferred. All scripts have been implemented in MATLAB 2019b and/or Mathematica 12 and have been run on an Intel\textregistered \ Core\texttrademark\ i7-4770 CPU 3.40GHz computer.

We shortly denote the PDE methods as RK4 and RK3b, where RK4 denotes the classical four-stage Runge-Kutta method \eqref{eq:RK4} combined with the first-order upwind space discretization, and RK3b the three-stage Runge-Kutta method \eqref{eq:second upwind space compact} combined with the third order accurate flux limited space discretization, as described in Section \ref{sec:2nd order upwind with flux}.

\subsection*{Example 1: \boldmath $S=2$}

In our first example we take a two-state Markov chain with generator matrix
\[\bm{Q}=\left(\begin{array}{rr}
-1 & 1\\
1 & -1
\end{array}\right)
.\]
The fluid rates are respectively $r_1 = 1$ and $r_2 = -2$ and the buffer capacity is $B^{\text{max}} = 10$. For this simple example, the LST has a closed-form expression:
\begin{align*}
\tilde{J}(x,w) &=
\frac{\left(\frac{w}{2}+1\right) e^{\frac{1}{4} x (-r(w)+w+1)} \left((r(w)-3 w-1) e^{\frac{1}{2} x r(w)}+e^{5 r(w)} (r(w)+3
   w+1)\right)}{(w+1) \left(e^{5 r(w)}+1\right) r(w)-\left(3 w^2+6 w+1\right) \left(1-e^{5 r(w)}\right)},
\end{align*}
with $r(w)=\sqrt{9 w^2+18 w+1}$. Still, the inversion of the LST has to be done numerically. We have compared the algorithms by Abate and Whitt \cite{AbateWhitt} and Den Iseger \cite{deniseger}, and in general the algorithm by Abate and Whitt is slightly faster, but it is also more susceptible to numerical errors for small values of $t$. It turns out that Den Iseger's algorithm is better at handling the jumps in the distribution function in $K_\alpha(x,t)$ (see also the discussion in Section \ref{sec:2nd order upwind with flux}), causing these numerical errors. It is straightforward to show that the smallest possible hitting time is $t=-x/r_S$ (since we assume that the rates are ordered), which equals $x/2$ in this example. Since we start in steady state, the probability mass at $t=x/2$ is equal to
\[
\pi_S e^{-q_{S,S}x/r_S}=\frac12 e^{-x/2}.
\]
In the sequel, we take $x=5$, resulting in a mean first passage time of $E[T]=11.27$. It is noteworthy that the LST inversion method does not have parameters $\Delta x$ or $\Delta t$, because it is possible to compute $K_\alpha(x,t)$ separately for each desired ($x, t$) combination. In contrast, the PDE methods (RK3b and RK4) compute $K_\alpha(x,t)$ on a large grid of ($x, t$) values and the accuracy of the algorithm depends on the distance between the grid points. For this reason we define two sets of parameters, which we refer to as ``low resolution'' and ``high resolution'':
\begin{center}
\begin{tabular}{|l|ll|}
\hline
Parameters & Low resolution & High resolution \\
\hline
$\Delta x$ & 0.1 & 0.01 \\
$\Delta t$ & 0.05 & 0.005 \\
$\Delta t$ LST & 0.1 & 0.01 \\
\hline
\end{tabular}
\end{center}
Note that we only choose $\Delta x$. We compute $\Delta t$ as $\Delta t =\displaystyle\min_{\alpha} \left|{\Delta x}/{r_\alpha}\right|$ for RK4, to ensure that the CFL condition is met, and $\Delta t =\frac14\displaystyle\min_{\alpha} \left|{\Delta x}/{r_\alpha}\right|$ for RK3b, as discussed below Equation \eqref{cfl}.
For the LST method we only need to specify $\Delta t$ to determine the points where we want to obtain the distribution function. In this example we use Abate and Whitt's algorithm to invert the LST. We evaluate the distribution function of $T$ for values of $t$ up to $t=20$.

\begin{figure}[ht]
\parbox{0.45\textwidth}{
\centering
\includegraphics[width=0.45\textwidth]{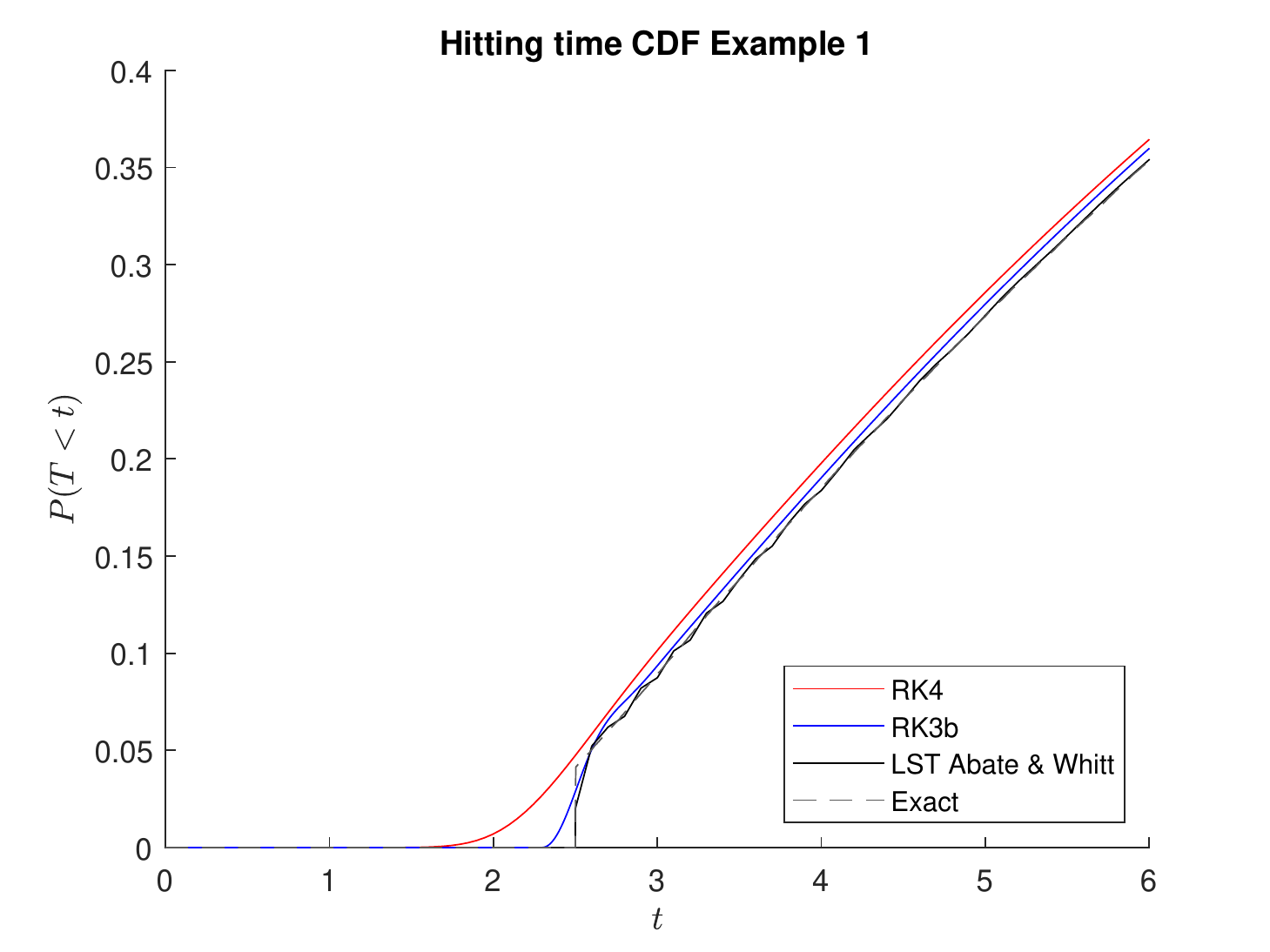} \\
(a) Low resolution}
\parbox{0.45\textwidth}{
\centering
\includegraphics[width=0.45\textwidth]{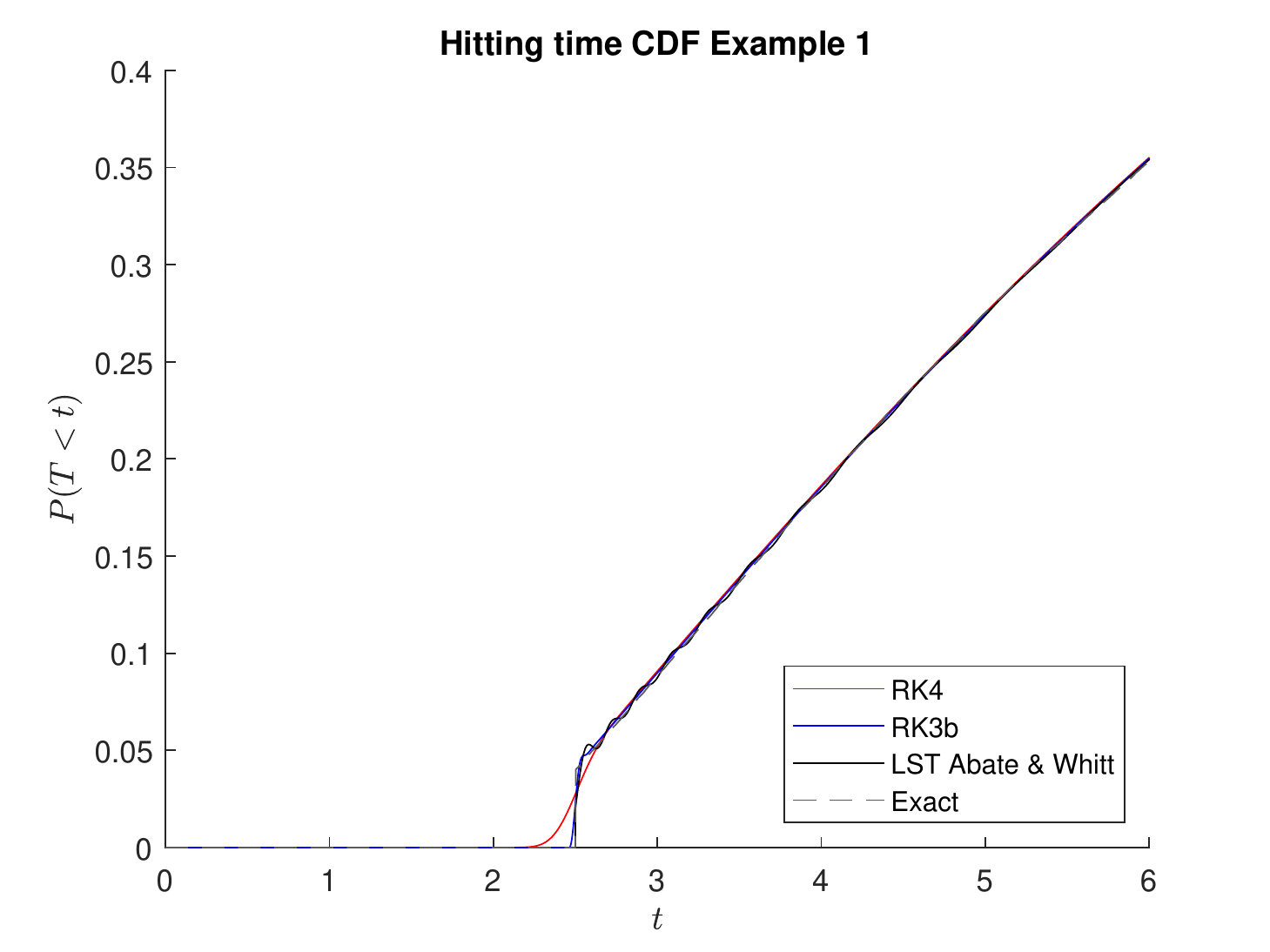}\\
(b) High resolution
}
\caption{$J(x,t)$ for Example 1, taking $x=5$.}
\label{fig:ex1}
\end{figure}

In Figure \ref{fig:ex1} we show how the accuracy of the algorithms is determined by the parameters, by comparing the numerical results with exact results obtained by simulation.
All methods suffer from loss of accuracy near the discontinuity at $t=2.5$.  The PDE methods perform very well, but we need a fine grid to obtain accurate results. Regarding the computational speed, it is clear that the RK4 method is by far the fastest method, but it is also the least accurate (with these parameter values). When one is interested in the distribution of the first passage time for one particular initial fluid level $x$ and the number of states in the background CTMC is limited to two, the LST method is clearly the optimal choice due to the existence of a closed-form expression for the LST.
In the next example we will consider a CTMC with more states and show that the oscillations near the discontinuity can be solved by using a different LST inversion method.

\begin{table}[ht]
\caption{Computing times (seconds) for Example 1.}
\label{tbl:compEx1}
\begin{center}
\begin{tabular}{|l|rrrr|}
\hline
Resolution & RK3b & RK4 & LST Symbolic & LST Numeric  \\
\hline
Low  & 1.10   & 0.02  &  0.09    &0.40     \\
High  & 45.41 &   0.28    &0.83    &3.17   \\
\hline
\end{tabular}
\end{center}
\end{table}

\subsection*{Example 2: \boldmath $S=5$}

\begin{figure}[ht]
\centering
\includegraphics[width=0.85\textwidth]{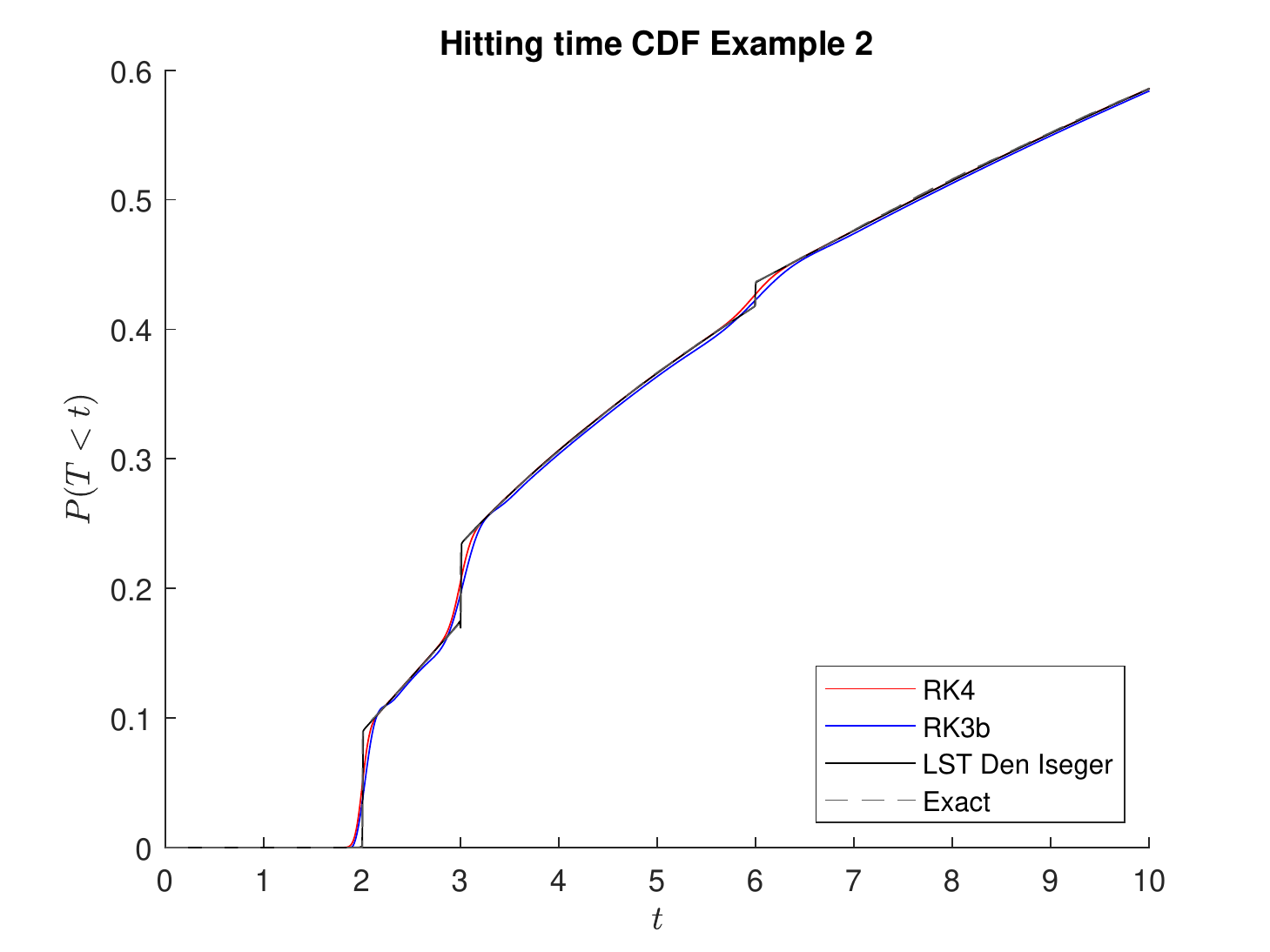} \\
\caption{$J(x,t)$ for Example 2, taking $x=6$.}
\label{fig:ex2}
\end{figure}

We use this example to compare the different methods when the LST is not of closed-form and to focus more on the discontinuities we observed in the previous example. Our background process now has five states, with generator matrix
\[\bm{Q}=\frac{1}{10}\left(\begin{array}{rrrrr}
-4 & 1 & 1 & 1 & 1\\
1 & -4 & 1 & 1 & 1\\
1 & 1 & -4 & 1 & 1\\
1 & 1 & 1 & -4 & 1\\
1 & 1 & 1 & 1 & -4\\
\end{array}\right)
.\]
We take the rates equal to $(4, 2, -1, -2, 3)$, the maximum buffer size $B^\text{max}=10$, and initial fluid level $x=6$. With these parameter values, the discontinuities are located at $t=2$, $t=3$, and $t=6$.
The state space of the CTMC is still relatively small ($S=5$), but already too large to allow for a simple, closed-form LST. As a consequence, we will only consider numerical LST inversion. We now only employ Den Iseger's algorithm, which is known to handle discontinuities better. Our main goal is to find the distribution function of $T$ for $t \leq 20$.
In contrast to the previous example, we now take parameter values that result in small differences between the computing times, which makes it interesting to compare the achieved accuracy of the different methods. The parameter values and the computing times are shown in Table~\ref{tbl:compEx2}. A plot of the computed cumulative distribution function of $T$ is shown in Figure~\ref{fig:ex2}. Note that this plot does not contain the Abate and Whitt inversion. Although this method outperforms Den Iseger in terms of speed, the oscillations near the discontinuities are so severe that we do not consider it to be a good alternative here. The PDE based methods also fail to perform well near these discontinuities, although they are still better than Abate and Whitt. Den Iseger is the only method that can really accurately approximate the true probability distribution. Due to the flexibility with respect to the $t$-values of the LST inversion methods, it is even possible to significantly decrease the computing time of the method by Den Iseger, by creating a ``flexible grid'' where the distribution function should be evaluated, since it really makes no sense (for most practical purposes) to use $\Delta t=0.01$ except near the discontinuities. We were able to reproduce Figure \ref{fig:ex2} by taking $\Delta t=0.01$ for $t < 7$ and $\Delta t=0.1$ for $t\geq 7$, reducing the computing time to only 3.31 seconds.

As a final disclaimer, we emphasize that the choice of parameter values in this example is completely motivated by the goal to get comparable computing times.
A more logical choice for the RK3b parameter values would be to take the same $\Delta t$ as for the RK4 PDE method ($\Delta x=0.02$ and $\Delta t=0.0012$). In that case, the accuracy would increase significantly, but the computing time would also increase, to 61 seconds.

\begin{table}[ht]
\caption{Parameters and computing times for Example 2. Abbreviations: LST AW = Abate and Whitt LST inversion, LST DI = Den Iseger LST inversion (with and without flexible grid).}
\label{tbl:compEx2}
\begin{center}
\begin{tabular}{|l|lllll|}
\hline
 & RK3b & RK4 & LST AW & LST DI & LST DI Flex Grid \\
\hline
$\Delta x$ & 0.1 & 0.005 &&&\\
$\Delta t$ & 0.0062 & 0.0012 & 0.01 & 0.01 & 0.01 / 0.1\\
Computing time (sec)  & 4.84&    5.80&    4.70&    7.32 & 3.31    \\
\hline
\end{tabular}
\end{center}
\end{table}

\subsection*{Example 3: \boldmath $S=100$ birth-death process}

\begin{figure}[ht]
\parbox{0.48\textwidth}{
\centering
\includegraphics[width=\linewidth]{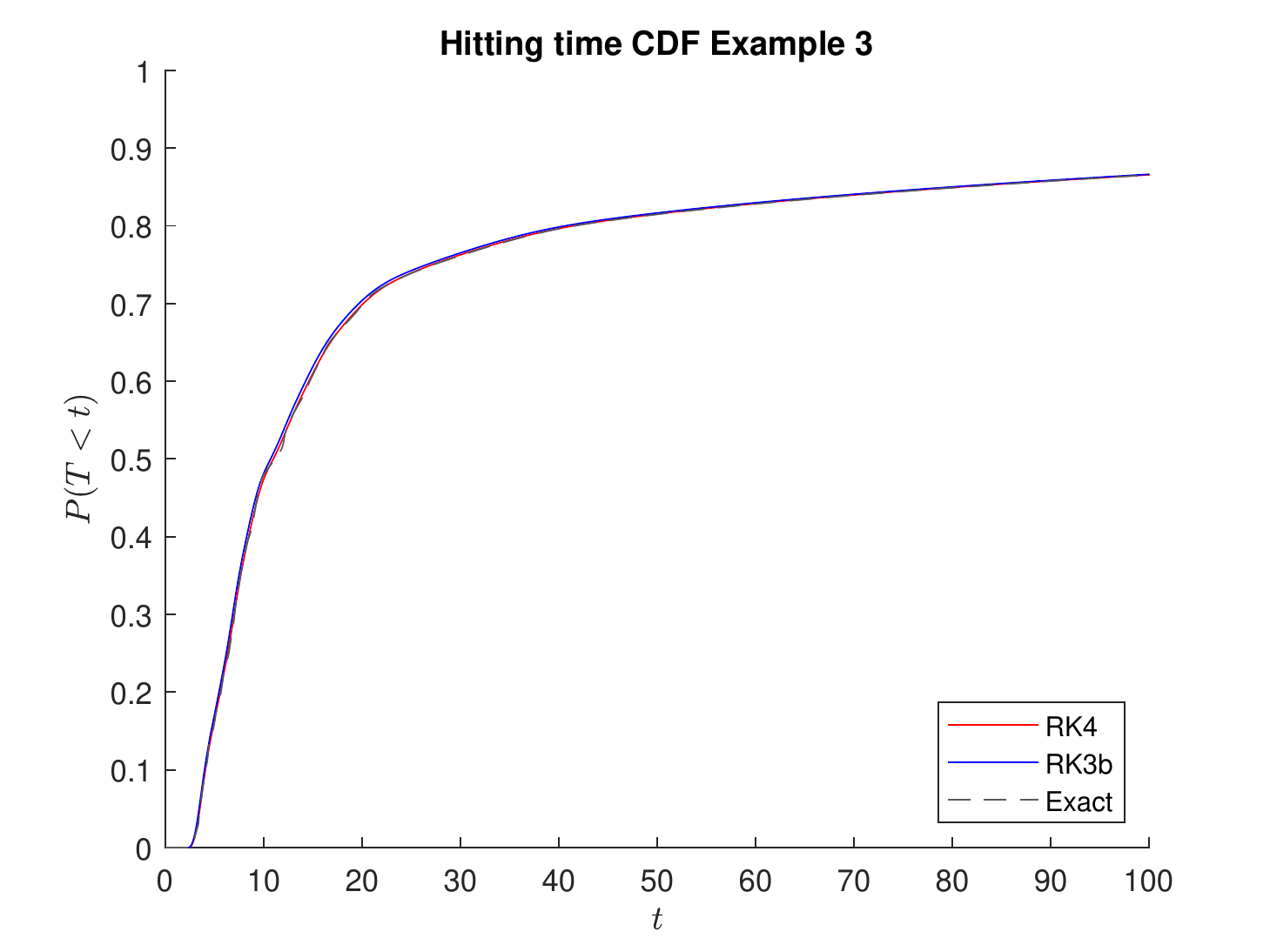} \\
(a) $J(x,t)$ }
\parbox{0.48\textwidth}{
\centering
\includegraphics[width=\linewidth]{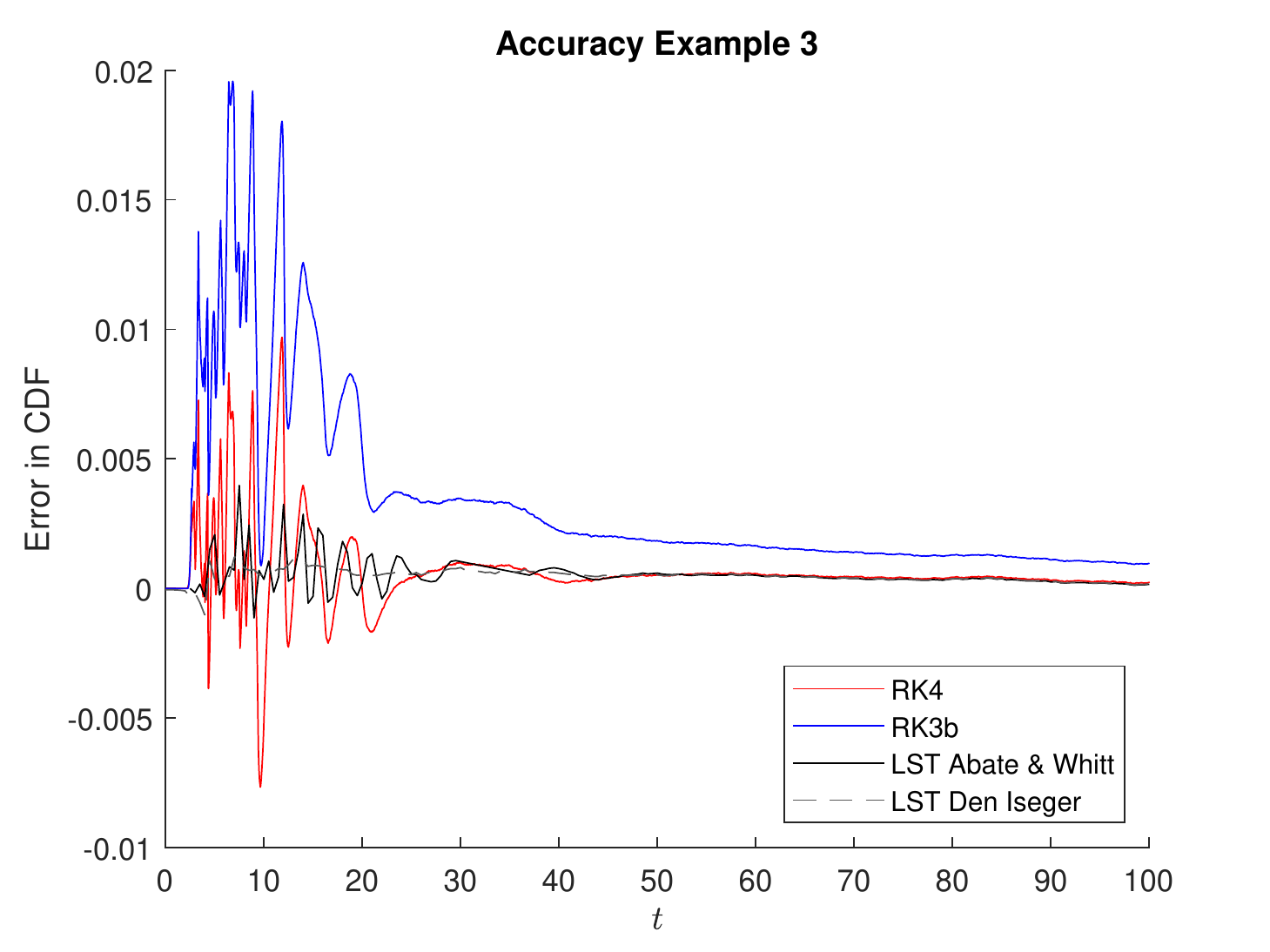} \\
(b) Error for $J(x,t)$ }
\caption{Numerical results for Example 3, taking $x=5$.}
\label{fig:ex3}
\end{figure}

In some applications, the underlying CTMC has a large number of states. For example, in \cite{debaratipaper,debaratiPhDthesis} fluid queues are used to model wind farm power output and their underlying Markov Chains have more than 100 states. In these settings, the discontinuities that we have been focussing on in the previous two examples, are less relevant from a practical perspective.
In this example, we generate a $\bm{Q}$ matrix for a birth-death process with $q_{\alpha\alpha}=-1$, $\forall\alpha\in\mathcal{S}$, $q_{\alpha,\alpha\pm1}=0.5$, $\forall\alpha\in\mathcal{S}\backslash\{1,100\}$ and periodic boundary.
The rates $r_\alpha$ are randomly  sampled from a normal distribution, i.e., $r_\alpha\sim\mathcal{N}(\mu,\sigma)$ with $\mu=-50$ and $\sigma=100$. We make sure that rates are sampled such that the drift is negative, i.e., $d < 0$. We take $B^\text{max} = 10$. First we compare $J(5, t)$ calculated from the PDE problem solver with the LST inversion methods and a Monte Carlo simulation.
As suggested before, for calculating $\Delta t$, we take
\[
\text{RK3b: }~\Delta t =\frac14\displaystyle\min_{\alpha\in\mathcal{S}} \left|\frac{\Delta x}{r_\alpha}\right|,   	
\qquad
\text{RK4: }~\Delta t =\displaystyle\min_{\alpha\in\mathcal{S}} \left|\frac{\Delta x}{r_\alpha}\right|,   	
\]
where the factor $\frac14$ for RK3b  is chosen to be more conservative than the CFL condition. 
With $\Delta x=0.05$ for RK4 and $\Delta x=0.5$ for RK3b, we find $\Delta t$ is equal to $0.0033$ and $0.033$, respectively. With these settings, the computing times for the PDE methods to determine the CDF for $0\leq t\leq 100$ are 52 seconds (RK3b) and  47 seconds (RK4). To assess the accuracy of the methods, we have obtained extremely accurate results from a Monte Carlo simulation, which we will refer to as ``exact'' in the sequel. Figure \ref{fig:ex3}(a) shows the CDF obtained by the two PDE methods and the Monte Carlo simulation.
It can be seen, in Figure \ref{fig:ex3}, that the accuracy is high, in particular for $t>20$. For smaller values of $t$ the PDE methods suffer from the same issues as discussed in the previous two examples. For most practical applications, however, one might argue that these inaccuracies are not of the utmost importance.

The LST method has not been considered in Figure \ref{fig:ex3}(a), because its accuracy is so high that it is practically indistinguishable from the exact values. However, it turns out to have quite some issues, partly caused by the large state space. First, there are numerical issues. The roots $s_k(w)$ vary greatly in size, which is magnified drastically by taking the exponential in \eqref{eq:til_H final}. This makes the numerical solution of \eqref{eq:systemofeqnsmatrixform} extremely unstable, to the extent that MATLAB cannot be used anymore and we have to use Mathematica instead. The second issue, which is partly related to the first issue, is that solving \eqref{eq:systemofeqnsmatrixform} becomes very slow when working with Mathematica's high-precision arithmetic. For $S=100$ the time to evaluate \emph{one single value} of the LST $\tilde{K}{_\alpha}(x,w)$ was approximately one second on the aforementioned Intel\textregistered \ Core\texttrademark\ i7-4770 CPU 3.40GHz computer. Considering that, with default parameter settings, Abate and Whitt's algorithm requires 28 evaluations of the LST to invert it at one single point, and Den Iseger's algorithm 64 LST evaluations, it is apparent that the PDE method tremendously outperforms the LST method with respect to the computing time.
The
computing times to obtain $J(x,t)$ for a single value of $x$, now with $\Delta t=0.5$, are approximately 1200 seconds for Abate and Whitt and 2900 seconds for Den Iseger, on a computer with multiple cores running the computations in parallel. Indeed, an advantage of the LST method is that it is extremely suitable for parallelization (we note that there may also be possibilities of parallelizing the PDE method; we leave exploration of such possibilities for future study).

Another advantage of the LST method is its accuracy, in particular when using Den Iseger's inversion algorithm. This is clearly confirmed by Figure \ref{fig:ex3}(b), which shows the errors of the various methods, i.e., the difference between the computed values and the exact values.
This example highlights the (dis)advantages of the different methods. The PDE methods are much faster than the LST method, but it is not straightforward to determine the right parameter values (the space/time grid) on beforehand. They lose accuracy where the CDF has discontinuities.
The accuracy can be improved by using a finer grid (i.e., smaller $\Delta x$ and $\Delta t$); this will increase the computing time.
With the LST method, in particular in combination with Den Iseger's inversion algorithm, one does not have to worry about parameter values and accuracy, provided that the numerical software to perform the calculations supports the required high-precision arithmetic.

\subsection*{Example 4: \boldmath $S=100$, non-sparse $Q$}

\begin{figure}[ht]
\parbox{0.48\textwidth}{
\centering
\includegraphics[width=\linewidth]{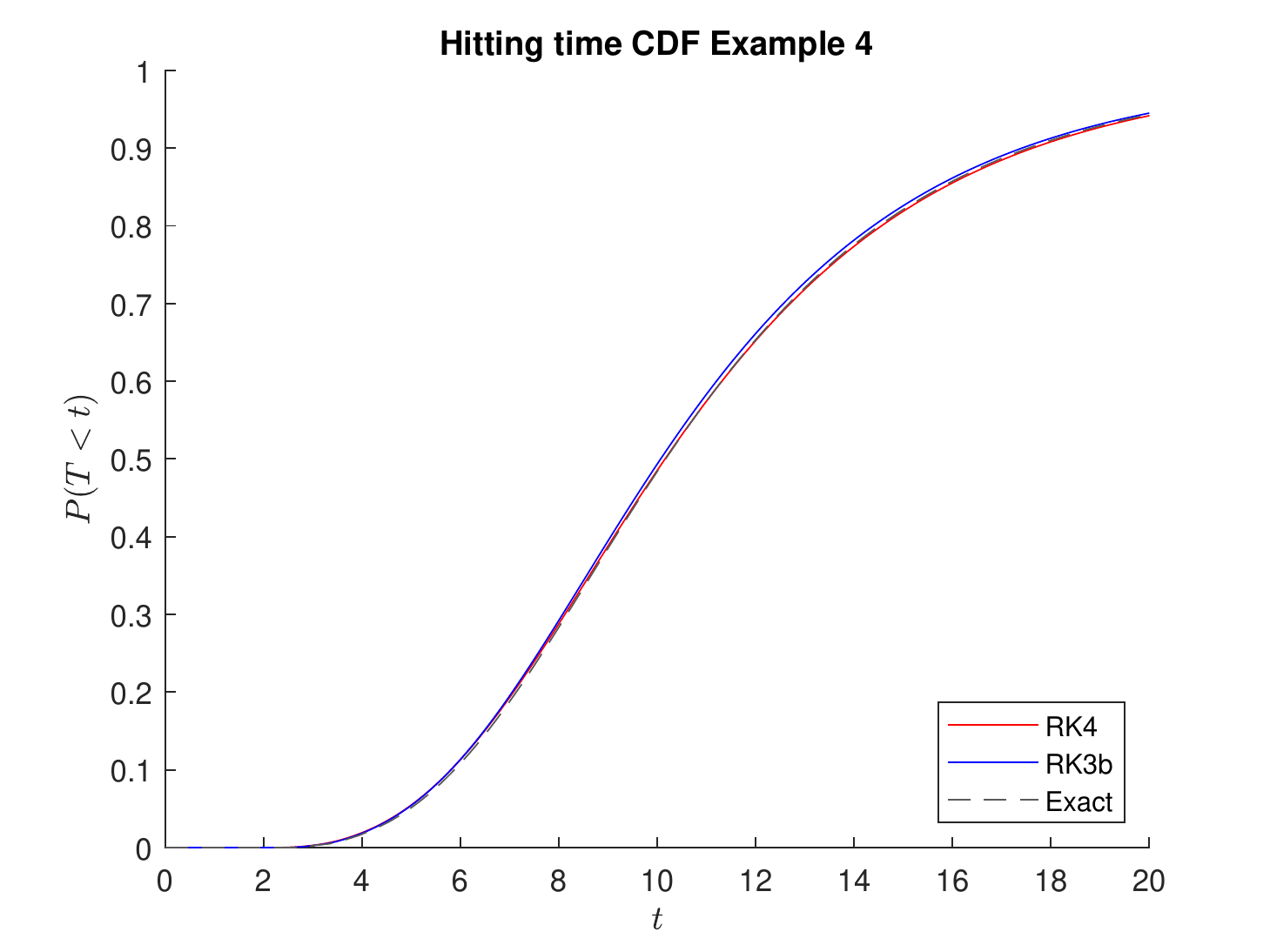} \\
(a) $J(x,t)$ }
\parbox{0.48\textwidth}{
\centering
\includegraphics[width=\linewidth]{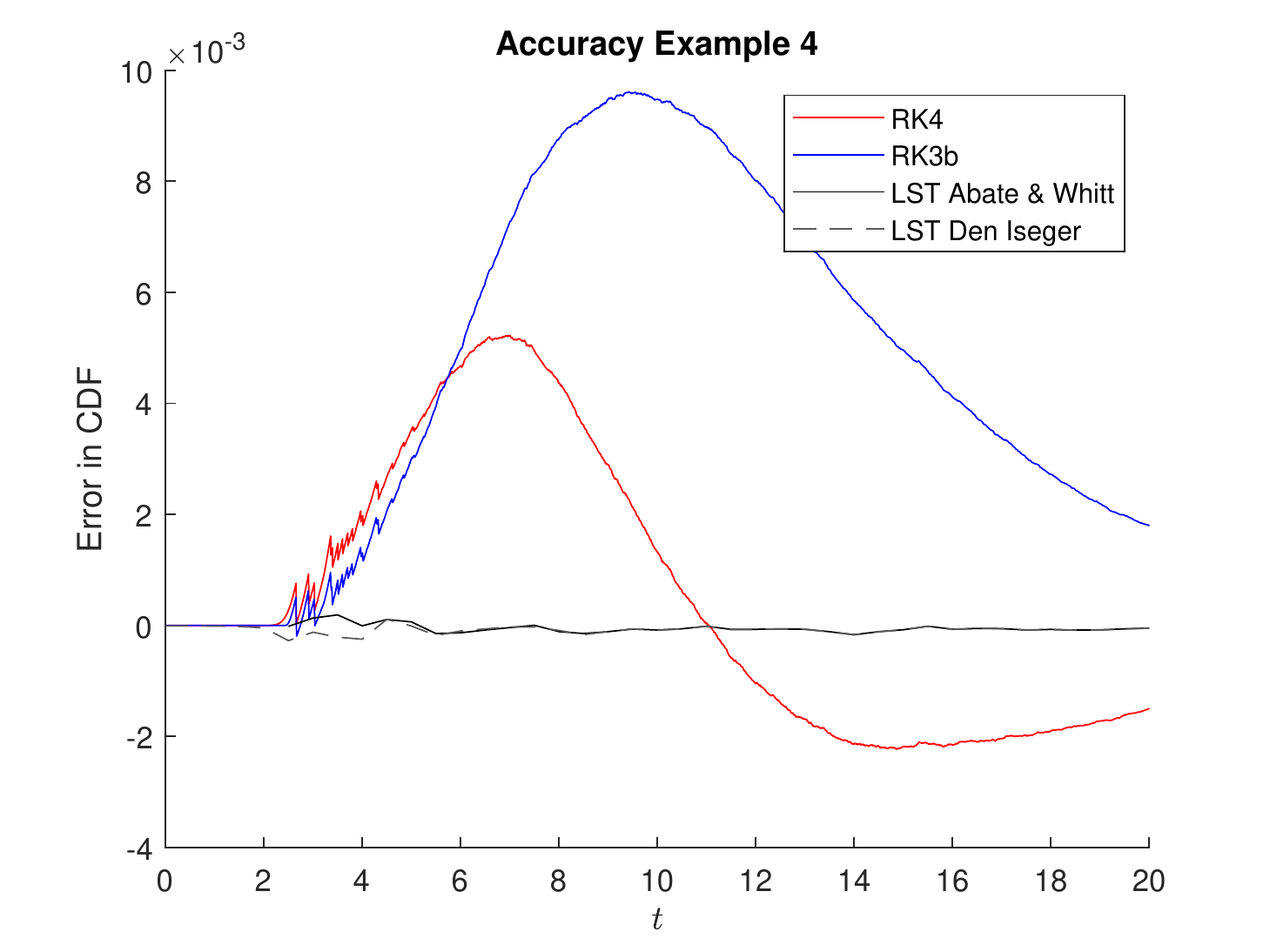} \\
(b) Error for $J(x,t)$ }
\caption{Numerical results for Example 4, taking $x=5$.}
\label{fig:ex4}
\end{figure}

In the previous example we used a $\bm{Q}$ matrix for a
birth-death process which is sparse in nature. In this example, we present results for a non-sparse $\bm{Q}$ matrix. More specifically, we take
\[
q_{\alpha\beta} = \begin{cases}
\frac{1}{100} & \qquad \alpha\neq\beta,\\
-\frac{99}{100} & \qquad \alpha=\beta.
\end{cases}
\]
We take the same rates as in the previous example. In order to reach a high accuracy and keep the computing times for RK3b and RK4 (more or less) similar, we take $\Delta t=0.0133$ for RK3b and $\Delta t=0.0033$ for RK4. Computing times are respectively 119 and 190 seconds, for RK3b and RK4. We take $\Delta t=0.5$ for the LST method to produce Figure \ref{fig:ex4}. The computing times for Abate \& Whitt and Den Iseger are 280 and 670 seconds, respectively. Again, we use parallel computations in Mathematica.

For this example, we can draw the same conclusions as in the previous examples. The main reason to include this example, is to have a setting where the discontinuities are so small that they are hardly noticeable. In this case, all methods perform relatively well. Obviously, the PDE methods are still faster than the LST methods, but since the CDF is so smooth, there is no reason to take a very fine grid. Moreover, the method by Abate and Whitt performs very well in this case.

\section{Conclusions}
\label{sec:conclusions}

Markov modulated fluid models are very relevant stochastic processes from both theoretical and practical point of view. As a consequence, they are well-studied by the stochastic community. However, determining the distribution of the first passage time turns out to be a computationally expensive numerical procedure. Indeed, this distribution function follows a hyperbolic partial differential equation, which is traditionally solved using methods based on Laplace-Stieltjes transforms, which need to be numerically inverted to obtain the desired probabilities. In this paper, we show that the PDE can also be solved directly using numerical integration methods. In particular, we use a third order accurate Runge-Kutta scheme in combination with a flux-limited third order accurate space discretization (see \cite{hundsdorfer1995positive}) and a fourth order accurate Runge-Kutta scheme in combination with the first order accurate upwind space discretization. These PDE-based methods have several advantages compared to the LST-based methods. In particular, the computing times are significantly smaller for Markov chains where the LST has no compact closed-form expression. Another major advantage of the PDE-based methods is that they compute the distribution of the first passage time for a very fine grid of values between zero and $B^{\text{max}}$, which saves even more computing time when one is interested in the first passage time for multiple initial buffer levels. The last advantage we discuss here, is that the PDE methods do not suffer from numerical issues that may arise with the LST method when solving the system of equations \eqref{eq:systemofeqnsmatrixform}.

There are also some disadvantages to the PDE methods. First, it is not clear on beforehand how fine the space-time grid should be to obtain accurate results. Another disadvantage is that $\Delta t$ and $\Delta x$ cannot be chosen independently. Even if one is interested in the first passage time for only a few $x$ values, $\Delta x$ should be chosen sufficiently small, for accuracy reasons. Finally, the PDE methods produce smooth curves, meaning that the accuracy near discontinuities in the first passage time distribution function is worse than in areas where the function is smooth. For most practical purposes, however, the PDE-based methods can give excellent results in significantly lower computing times.

As a final remark, we note that development of numerical PDE solvers is a large and mature field of research. Over the years, many methods, variations and refinements have been proposed. We have only considered a few well-established and rather straightforward methods in this study. We anticipate that there are various possibilities of further improving and refining the use of PDE methods for computing first passage times. This can lead to more accurate and/or more computationally efficient algorithms, as well as better handling of discontinuities and more automated error control. Some examples are non-uniform grids for spatial discretization, mesh refinement, adaptive timestepping and parallelization. We leave exploration of such more sophisticated methods in the context of first passage time problems for future study.

\section*{Acknowledgments}

This work was part of the Industrial Partnership Program (IPP) ``Computational Sciences for Energy Research'' of the Foundation for Fundamental Research on Matter (FOM), which is now part of the Netherlands Organization for Scientific Research (NWO). This research program was co-financed by Shell Global Solutions International B.V..

\bibliographystyle{abbrv}

\appendix
\section*{Appendix}
\section{First order accurate upwind scheme}\label{App:first order upwind}

The matrix $\bm{A}$ in \eqref{eq:upwind space compact} can be written as $\bm{A} = \tilde{\bm{R}} + \tilde{\bm{Q}} + \bar{\bm{Q}}$, in which $\tilde{\bm{R}}$ is
a block diagonal matrix given by $\tilde{\bm{R}} = \verb|blkdiag|(\tilde{\bm{R}}_1,\dots, \tilde{\bm{R}}_S)$, where

\begin{align*}
\tilde{\bm{R}}_\alpha &=
\frac{r_\alpha}{\Delta x}
\left[
\begin{array}{rrrrr}
0&0&0&0& \hdots \\
-1& 1 & 0&0&\hdots \\
0&-1&1&0&\dots\\
\vdots&\vdots&\vdots&\vdots&\ddots\\
0&0&0& \hdots &0\\
\end{array}
\right]
\qquad &&\forall\{\alpha : r_\alpha < 0\} \qquad \text{and, }\\
\tilde{\bm{R}}_\alpha &=
\frac{r_\alpha}{\Delta x}
\left[
\begin{array}{rrrrr}
0&0&0&0& \hdots \\
0&-1&1&0&\dots\\
\vdots&\vdots&\vdots&\vdots&\ddots\\
0& 0 & 0&\hdots &0\\
\end{array}
\right]
\qquad &&\forall\{\alpha : r_\alpha > 0\}.
\end{align*}

Furthermore, $\tilde{\bm{Q}}={\bm{Q}} \otimes {\bm{I}}_{q\times q}$ where ${\bm{I}}_{q\times q}$ is the identity matrix of size $q \times q$. The matrix $\bar{\bm{Q}}$
contains the discretization of the boundary conditions
for $K_\alpha^{n,q}$ for $\{\alpha : r_{\alpha}>0\}$
and has dimensions equal to $\tilde{\bm{Q}}$. The
elements of the matrix are given by
$\bar{Q}_{n \alpha,n\beta}=\tilde{Q}_{n \alpha,n\beta}$ for  $\{\alpha : r_{\alpha}>0\}$
and $\beta=1,\dots,S$. All the other elements of matrix $\bar{\bm{Q}}$
are zero.

\section{Flux limited higher order accurate upwind scheme}
\label{App:second order upwind}

We present the details of the flux limited higher order accurate upwind scheme here.
From \eqref{eq:flux limiter} it can be seen that flux values for faces $F_\alpha^{3/2}$
and $F_\alpha^{n+1/2}$, $\forall \{\alpha : r_\alpha \leq 0\}$  cannot be defined by \eqref{eq:flux limiter}. Similarly, flux values for faces $F_\alpha^{1/2}$
and $F_\alpha^{n-1/2}$, $\forall \{\alpha : r_\alpha > 0\}$
cannot be defined by \eqref{eq:flux limiter}. Hence we define
the fluxes for all the faces as follows:

\begin{itemize}
  \item $\bm{\forall \{\alpha : r_\alpha \leq 0\}}$
  \\
  \begin{enumerate}
  \item $F_\alpha^{1/2}$ is given by the boundary condition in Section \ref{sec:first passage time} .
  \item $F_\alpha^{3/2} = \frac12(K^1_\alpha+K^2_\alpha)$, i.e., the central scheme.
  \item $F_\alpha^{p+1/2}$ for $p=2,\dots,n-1$ is given by \eqref{eq:flux limiter}.
  \item $F_\alpha^{n+1/2}= (\frac32K^n_\alpha-\frac12K^{n-1}_\alpha))$, i.e., the second order accurate upwind scheme
without limiters.
\end{enumerate}
  \item $\bm{\forall \{\alpha : r_\alpha > 0\}}$
  \\
  \begin{enumerate}
  \item $F_\alpha^{1/2}=(\frac32K^1_\alpha-\frac12K^{2}_\alpha))$ .
  \item $F_\alpha^{p+1/2}$ for $p=1,\dots,n-2$ is given by \eqref{eq:flux limiter}.
  \item $F_\alpha^{n-1/2} = \frac12(K^n_\alpha+K^{n-1}_\alpha)$, i.e., the central second order accurate scheme.
  \item $F_\alpha^{n+1/2}$  is given by the boundary condition in Section \ref{sec:first passage time}.
\end{enumerate}
\end{itemize}

Note that in \eqref{eq:flux limiter}, the expression $\phi_{p+1/2}=\phi(f_{p+1/2})$ was still left unspecified. As a convenience we repeat here how we specified them later in
Section \ref{sec:flux limiters}. As limiter $\phi$ we chose the Koren limiter \eqref{eq:barry limiter}, yielding third order accuracy, and $f_{p+1/2}$ is defined as
\begin{equation}
f_{p+1/2}=\begin{cases}
\frac{K_\alpha^{p+1}-K_\alpha^{p}}{K_\alpha^{p}-K_\alpha^{p-1}} &\qquad  \forall \{\alpha : r_\alpha \leq 0\},   \\
\frac{K_\alpha^{p-1}-K_\alpha^{p}}{K_\alpha^{p}-K_\alpha^{p+1}} &\qquad  \forall \{\alpha : r_\alpha > 0\}.   \\
\end{cases}
\end{equation}

\end{document}